\documentclass[a4paper,conference]{IEEEtran}
\IEEEoverridecommandlockouts

\usepackage{cite,color,url}
\usepackage{amsmath,amssymb,amsfonts,hyperref,multirow}
\usepackage{booktabs}
\usepackage{algorithmic,tikz}
\usepackage{subcaption,amsthm,longtable,lscape,amsmath}
\usepackage{graphicx,longtable}
\usepackage{textcomp}
\usepackage{xcolor}
\usetikzlibrary{fit,calc}
\def\BibTeX{{\rm B\kern-.05em{\sc i\kern-.025em b}\kern-.08em
    T\kern-.1667em\lower.7ex\hbox{E}\kern-.125emX}}
\begin{document}

\title{On Solving the Assignment Problem with Conflicts}
\author{\IEEEauthorblockN{Roberto Montemanni}
\IEEEauthorblockA{\textit{Department of Sciences and Methods for Engineering} \\
\textit{University of Modena and Reggio Emilia}\\
Reggio Emilia, Italy \\
roberto.montemanni@unimore.it}
\and
\IEEEauthorblockN{Derek H. Smith}
\IEEEauthorblockA{\textit{Faculty of Computing, Engineering and Science} \\
\textit{University of South Wales}\\
Pontypridd, Wales, UK \\
derek.smith@southwales.ac.uk}}

\maketitle
\begin{abstract}
A variant of the well-known Assignment Problem is studied in this paper, where pairs of assignments are conflicting, and cannot be selected at the same time. This configures a set of hard constraints. The problem, which models real applications, looks for a complete assignment that minimizes the total cost, while no conflict is violated. 

In this paper, we consider a previously known mixed integer linear
program representing the problem and we solve it with the open-source solver CP-SAT, part of the Google OR-Tools computational suite.

An experimental campaign on the instances available from the literature, indicates that the approach we propose achieves results comparable with, those of state-of-the-art solvers, notwithstanding its intrinsic conceptual and implementation simplicity. The solver adopted is also able to provide heuristic solutions quicker and better than the heuristic methods previously discussed in the literature.
\end{abstract}


\section{Introduction} \label{into}
The Assignment Problem is a well-known optimization problem with direct applications in several fields such as  personnel scheduling, task
assignment, job shop loading, facility location and workforce planning. Instances of any practical size can  be solved efficiently in polynomial time
using, for example the Hungarian algorithm \cite{ahu93}, \cite{kuh55}.
 
The focus of this work is on the Assignment Problem with Conflicts (APC), a problem introduced recently in \cite{onc19}, which is an extension of AP with additional conflict
constraints that forbid pairs of assignments from happening at the same time, imposing therefore hard constraints. The APC deals with finding a minimum weight perfect matching such that no more than one edge is selected from
each conflicting edge pair. 

Classic optimization problems with conflicts have been studied extensively in the last decades, due to their practical implications, and their extended computational complexity in their general settings over the standard problems. For example, it is possible to trace studies on knapsack problems with conflicts \cite{bet17}, \cite{con21}, on spanning trees with conflicts \cite{zha11}, \cite{som14}, \cite{car19}, \cite{car19b}, \cite{car21}, on shortest paths with conflicts \cite{cer23}, set coverings with conflicts \cite{saf22}, \cite{car24} and maximum flows with conflicts \cite{suv20}, \cite{car25}. 

There are several practical applications of the APC. For example in container terminals, where arriving containers ($V_A$) have to be organized in the yard locations ($V_B$). Containers with certain characteristics -- typically related either to size, weight or priorities -- cannot be stacked on top of each others, leading therefore to conflict assignments. The costs of the assignment can be related to other real-world factors such as distance from the ship and the yard location. These settings configure an APC to be solved in order to optimize yard operations. Other similar applications of the APC are in general warehouses or multi-compartment vehicles, where goods with certain characteristics cannot be stored next to each other (e.g. food and toxic products). 
In personnel scheduling problems (e.g. airlines rostering) there might be incompatibilities between people, that should therefore not be assigned to the same team. This can be translated into straightforward conflicts within an APC.

Theoretical results on the computational complexity of the APC, together with some approximation results, can be found in \cite{dar11}. Further results on special polynomially-solvable settings were presented in \cite{onc13}, together with some heuristics and lower bounding schema. 
More contribution for the APC can be found in \cite{onc17} and \cite{onc18}, where some Mixed Integer Linear Programming models, exact Branch\&Bound methods and heuristic ideas are introduced. The concepts of the last two papers were later summarized and extended in \cite{onc19}, which remains the reference work for the problem under investigation.

In this paper, a mixed integer linear programming model for the APC is considered and solved by the open-source solver CP-SAT, which is part of the Google OR-Tools \cite{cpsat} optimization suite. Successful application of this solver on optimization problems with  characteristics similar to the problem under investigation, motivated our study \cite{md23}, \cite{cor}, \cite{rm25}. 
An experimental campaign on the  benchmark instances previously proposed in the recent literature is also presented and discussed.

The overall organization of the paper can be summarized as follows. The Assignment Problem with Conflicts is formally defined in Section \ref{desc}. Section \ref{model} discusses a mixed integer linear programming model to represent the problem. In Section \ref{exp} the approach we propose is compared with recent state-of-the-art methods from the literature. Conclusions are finally drawn in Section \ref{conc}.

\vspace{1cm}
\section{Problem Description}\label{desc}
The Assignment Problem with Conflicts (APC) can be formally described as follows.  Let $G = (V_A \cup V_B, E)$ be a complete bipartite graph, where $|V_A|=|V_B|$, $V_A \cap V_B = \emptyset$, and $E = \left \{ \{i,j\} | i \in V_A, j \in V_B \right \}$ is the set of possible assignments. A non-negative cost $c_{ij}$ is given for each $\{i,j\} \in E$. Moreover, a set $C= \left \{  \{i, j\}, \{k, l\} | \left \{ \{i, j\}, \{k, l\} \right \} \in E \right \}$ represents the conflicts among pairs of assignments. 

The objective of the APC is to find the perfect matching of minimum cost between the sets $V_A$ and $V_B$ that uses at most one edge from every conflict of $C$.

\begin{figure*}[h!]
{
\begin{center}
{
\begin{tikzpicture}[node distance={2cm}, main/.style = {draw, circle}]
			\node[main,minimum size=1.15cm] (0) {\textcolor{cyan}{$0$}};
			\node[main,minimum size=1.15cm] (1) [below of=0] {\textcolor{cyan}{$1$}};
			\node[main,minimum size=1.15cm] (2) [below of =1] {\textcolor{cyan}{$2$}};
			\node[main,minimum size=1.15cm] (3) [below of=2 ] {\textcolor{cyan}{$3$}};
			\node[main,minimum size= 1.15cm] (4) [below of=3] {\textcolor{cyan}{$4$}};
			\node[main,minimum size= 1.15cm] (10) [color=white,right of=0] {};
			\node[main,minimum size=1.15cm] (5) [right of=10]{\textcolor{magenta}{0}};
			\node[main,minimum size=1.15cm] (6) [below of=5] {\textcolor{magenta}{1}};
			\node[main,minimum size=1.15cm] (7) [below of =6] {\textcolor{magenta}{2}};
			\node[main,minimum size=1.15cm] (8) [below of=7] {\textcolor{magenta}{3}};
			\node[main,minimum size= 1.15cm] (9) [below of=8] {\textcolor{magenta}{4}};
			\draw [color=red, thick, line width=1.2,-] (0) to node {} (5) ;
			\draw [color=red, thick, line width=1.2,-] (4) to node {} (5) ;
			\draw [color=blue, thick, line width=1.2,-] (0) to node {} (7) ;
			\draw [color=blue, thick, line width=1.2,-] (3) to node {} (6) ;
			\draw [color=orange, thick, line width=1.2,-] (3) to node {} (5) ;
			\draw [color=orange, thick, line width=1.2,-] (4) to node {} (7) ;
			\draw [color=green, thick, line width=1.2,-] (1) to node {} (6) ;
			\draw [color=green, thick, line width=1.2,-] (2) to node {} (9) ;
\end{tikzpicture}
\hspace{2cm}
\begin{tikzpicture}[node distance={2cm}, main/.style = {draw, circle}]
			\node[main,minimum size=1.15cm] (0) {\textcolor{cyan}{$0$}};
			\node[main,minimum size=1.15cm] (1) [below of=0] {\textcolor{cyan}{$1$}};
			\node[main,minimum size=1.15cm] (2) [below of =1] {\textcolor{cyan}{$2$}};
			\node[main,minimum size=1.15cm] (3) [below of=2 ] {\textcolor{cyan}{$3$}};
			\node[main,minimum size= 1.15cm] (4) [below of=3] {\textcolor{cyan}{$4$}};
			\node[main,minimum size= 1.15cm] (10) [color=white,right of=0] {};
			\node[main,minimum size=1.15cm] (5) [right of=10]{\textcolor{magenta}{0}};
			\node[main,minimum size=1.15cm] (6) [below of=5] {\textcolor{magenta}{1}};
			\node[main,minimum size=1.15cm] (7) [below of =6] {\textcolor{magenta}{2}};
			\node[main,minimum size=1.15cm] (8) [below of=7] {\textcolor{magenta}{3}};
			\node[main,minimum size= 1.15cm] (9) [below of=8] {\textcolor{magenta}{4}};
			\draw [color=red, thick, line width=1.2,-] (4) to node {} (5) ;
			\draw [color=blue, thick, line width=1.2,-] (0) to node {} (7) ;
			\draw [color=green, thick, line width=1.2,-] (1) to node {} (6) ;
			\draw [color=black, thick, line width=1.2,-] (2) to node {} (8) ;
			\draw [color=black, thick, line width=1.2,-] (3) to node {} (9) ;
\end{tikzpicture}
}
	\caption{On the left an example of a APC instance is presented, where conflicts are indicated as colored edges. Costs are omitted. On the right, a feasible solution is represented. Observe that at most one edge is used for each color (conflict) and some black edges (not involved in any conflicts) are added to complete the matching.}
	\label{figu}
\end{center}
}
\end{figure*}
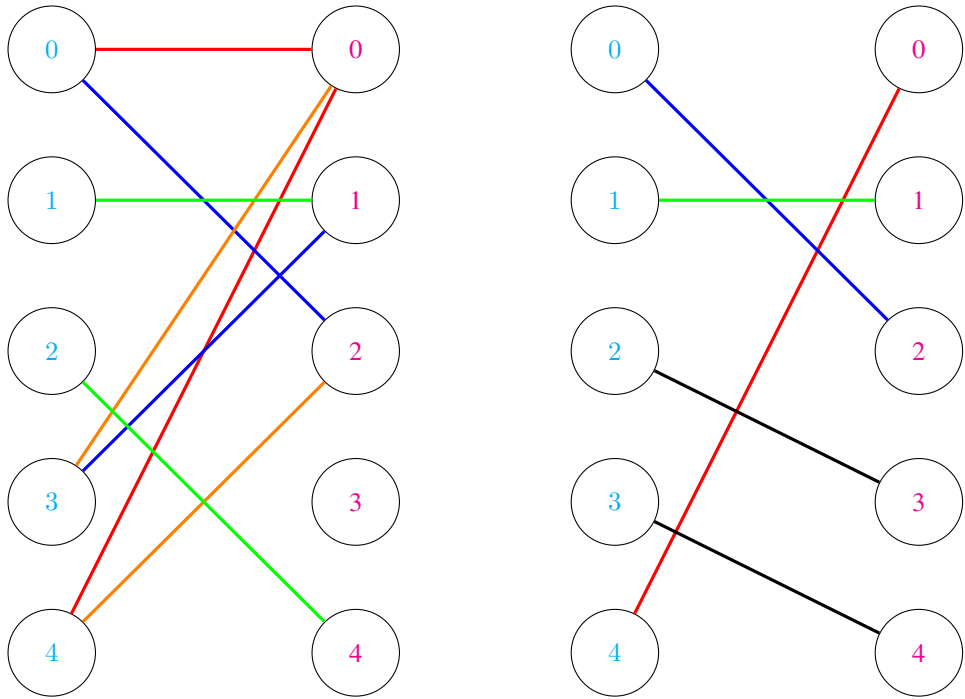

A small example of an APC instance and a respective feasible solution are depicted in Figure \ref{figu}.

\section{A Mixed Integer Linear Programming Model}\label{model}
In this section, a model for the APC, previously discussed in \cite{onc19}, is presented. A variable $x_{ij}$ takes value 1 if the edge between $ i\in V_A$ and $ j\in V_B$ is selected, 0 otherwise. The resulting model is as follows.

\begin{align} 
  \min \ \ &   \sum_{\{i,j\} \in E} c_{ij} x_{ij}& \label{1}\\ 
s.t. \ \ &	\sum_{j: \{i,j\} \in E}  x_{ij} = 1& i \in V_A \label{2}\\
	&	\sum_{i: \{i,j\} \in E}  x_{ij} = 1& j \in V_B \label{3}\\
	& x_{ij} + x_{kl} \le 1& \left \{ \{i, j\}, \{k, l\} \right \} \in E \label{4}\\
& x_{ij} \in \{0,1\} & \{i,j\} \in E \label{5}
\end{align}

The objective function (\ref{1}) minimizes the cost of the assignment selected.
Constraints (\ref{2}) impose that each element of $V_A$ has to be assigned to exactly one element of $V_B$ through the feasible edges. 
Constraints (\ref{3}) impose that each element of $V_B$ has to be assigned to exactly one element of $V_A$ through the feasible edges. 
Inequalities (\ref{4}) model the conflicts, imposing that at most one of two conflicting edges can be selected.
The domains of the $x$ variables are finally specified in constraints (\ref{5}).

\section{Computational Experiments} \label{exp}
In Section \ref{ben} we describe the benchmark instances previously introduced in the literature, and used for the present study. In Section \ref{res} the approach we propose is compared with the other methods available in the literature.

\subsection{Benchmark Instances}\label{ben}
In the literature, the only available benchmark set for the APC is -- to our knowledge -- the one proposed in \cite{onc19}. We will therefore adopt these instances for our experiments. 
The number of nodes in the left hand-side of the bipartite graph $G= (V1(G) \cup V2(G), E(G))$, indicated as $|V(G)|$, takes values between 15 and 500. Conflict pairs are generated at random,  and the number $|E(C)|$ of such pairs is between 5000 and 200000. A total of 135 test problems were introduced, but only 130 can be used for the experiments reported in this work, due to some inconsistencies in the available dataset. 

We refer the interested reader to \cite{onc19} for a comprehensive description of the instances.

\subsection{Experimental Results} \label{res}
The model discussed in Section \ref{model} has been solved with the Google OR-Tools CP-SAT solver \cite{cpsat} version 9.12. The experiments have been run on a computer equipped with an Intel Core i7 12700F CPU. The experiments for the methods previously appeared in \cite{con21} and against which we compare, were run on an a machine equipped with a 2.2 GHz Intel Core i7 processor (no more information is available), which according to \url{http://gene.disi.unitn.it/test/cpu_list.php} is  3 or more times slower.  

The methods involved in the comparison are:
\begin{itemize}
\item LS: Local Search heuristic algorithm discussed in \cite{onc19} (the original concepts of the method had already been discussed in \cite{onc18});
\item RDS: Russian Doll Search heuristic algorithm discussed in \cite{onc19};
\item BIP: best results obtained by solving with CPLEX 12.7 \cite{cplex} the two Binary Linear Programs presented in \cite{onc19} (the model had already been introduced in \cite{onc17} and \cite{onc18});
\item B\&B: Branch-and-Bound approach (with the best settings) presented in \cite{onc19} (a preliminary version of the method appeared in \cite{onc18});
\item CP-SAT: the mixed integer linear program presented in Section \ref{model} solved with Google OR-Tools CP-SAT solver 9.12 \cite{cpsat}.
\end{itemize}

The results are summarized in Table \ref{t1}, where for each group of five instances identified by the values of $|V_A|$ and $|C|$, the average of the optimal solutions is reported ($Opt$). All the methods considered were run for a maximum of 3600 seconds on each instance, and for each group of instances the following data are reported:
\begin{itemize}
\item Gap \%: the average optimality, calculated as $100~\cdot~\frac{Val-Opt}{Opt}$, where $Val$ is the value returned by the method. This value is reported only for heuristic methods;
\item Sec Best: the average computation time in seconds required to retrieve the best solution found. This value is reported only for heuristic methods;
\item Sec Tot: the average computation time in seconds required to prove the optimality of the best solution found. This value is reported only for exact methods.
\end{itemize}
The last lines of the table contain the averages of the indicators tabled for  the different methods.

\begin{table*}[ht]
{\begin{center}
\caption{Computational results.}\label{t1}
\begin{tabular}{lrrrrrrrrrr}
\toprule
\multicolumn{2}{l}{Instances}             & Opt     & \multicolumn{2}{c}{LS \cite{onc19}}          &  \multicolumn{2}{c}{RDS \cite{onc19}}      & BIP \cite{onc19}    & B\&B \cite{onc19}   & \multicolumn{2}{c}{CP-SAT}        \\
    $|V_A|$        &  $|C|$             &         & Gap \% & Sec Best & Gap \% & Sec  Best  & Sec  Opt  & Sec  Opt  & Sec Best & Sec Opt \\
\cmidrule(lr){1-2}\cmidrule(lr){3-3}\cmidrule(lr){4-5}\cmidrule(lr){6-7}\cmidrule(lr){8-8}\cmidrule(lr){9-9}\cmidrule(lr){10-11}
15          & 5000          & 2246.2  & 1.81   & 0.7  & 0.00   & 3.6    & 13.2   & 1.4    & 0.2      & 1.7     \\
20          & 10000         & 2450.8  & 1.99   & 0.8  & 0.00   & 13.1   & 130.1  & 6.3    & 17.5     & 886.7   \\
30          & 20000         & 3247.0    & 2.13   & 0.9  & 0.00   & 11.5   & 25.8   & 3.0    & 21.5     & 117.5   \\
30          & 30000         & 3355.8  & 2.18   & 1.1  & 0.00   & 637.4  & 1650.3 & 129.0  & 623.2   & 1931.1  \\
40          & 40000         & 4205.0    & 2.17   & 1.8  & 0.00   & 24.7   & 12.2   & 3.5    & 17.9     & 56.5    \\
50          & 50000         & 5183.4  & 2.10   & 2.5  & 0.00   & 39.4   & 7.6    & 3.3    & 24.4     & 38.8    \\
50          & 60000         & 5186.2  & 2.22   & 2.9  & 0.00   & 197.5  & 14.9   & 3.9    & 36.8     & 68.6    \\
60          & 80000         & 6163.6  & 2.25   & 3.7  & 0.00   & 77.6   & 10.6   & 5.6    & 42.2     & 60.2    \\
70          & 100000        & 7156.8  & 2.31   & 3.9  & 0.00   & 54.7   & 9.2    & 2.6    & 29.6     & 38.4    \\
70          & 150000        & 7166.8  & 2.44   & 4.1  & 0.01   & 2070.7 & 46.7   & 104.4  & 118.2    & 204.8   \\
80          & 200000        & 8154.2  & 2.43   & 4.3  & 0.05   & 2667.9 & 34.4   & 141.0  & 68.8     & 104.9   \\
90          & 250000        & 9157.2  & 2.45   & 4.3  & 0.07   & 3599.6 & 68.5   & 1126.9 & 223.1    & 359.5   \\
100         & 100000        & 10118.8 & 2.47   & 4.2  & 0.00   & 6.5    & 7.8    & 0.9    & 11.4     & 12.3    \\
100         & 250000        & 10144.0   & 2.49   & 4.5  & 0.01   & 1242.8 & 23.9   & 875.5  & 72.1     & 86.6    \\
100         & 350000        & 10160.0   & 2.61   & 4.7  & 0.42   & 3599.4 & 149.9  & 3171.5 & 451.0    & 832.6  \\
150         & 200000        & 15093.6 & 2.87   & 5    & 0.00   & 1.3    & 15.3   & 0.4    & 22.0     & 22.5    \\
150         & 350000        & 15102.2 & 2.93   & 5.1  & 0.00   & 117.5  & 22.7   & 7.1    & 32.9     & 34.1    \\
150         & 500000        & 15105.8 & 3.16   & 6.5  & 0.00   & 337.1  & 33.8   & 38.1   & 54.5     & 55.9    \\
200         & 200000        & 20077.4 & 3.54   & 8.2  & 0.00   & 0.4    & 39.3   & 0.1    & 25.5     & 26.5    \\
200         & 400000        & 20082.8 & 3.68   & 8.6  & 0.00   & 6.5    & 46.5   & 1.7    & 38.5     & 39.6    \\
250         & 500000        & 25058.2 & 4.08   & 10.5 & 0.00   & 6.6    & 99.5   & 0.9    & 50.2     & 51.4    \\
250         & 700000        & 25057.6 & 4.33   & 11.3 & 0.00   & 2.9    & 106.8  & 2.5    & 62.4     & 64.2    \\
300         & 100000        & 30042.4 & 4.49   & 12.8 & 0.00   & 0.2    & 184.1  & 0.1    & 31.9     & 32.9    \\
300         & 300000        & 30040.6 & 4.72   & 16   & 0.00   & 0.6    & 188.6  & 0.1    & 40.7     & 41.9    \\
400         & 200000        & 40016.6 & 5.16   & 20.5 & 0.00   & 0.3    & 608.8  & 0.0    & 57.9     & 59.7    \\
500         & 200000        & 50006.2 & 7.12   & 23.7 & 0.00   & 6.9    & 1497.3 & 0.0    & 85.7     & 87.7    \\
\midrule
\multicolumn{3}{l}{Averages} &          3.08   & 6.6  & 0.02   & 566.4  & 194.1  & 216.5  & 141.6    & 204.5  \\
\bottomrule
\end{tabular}
\end{center}}
\end{table*}

The results suggest the following considerations. First, the Local Search heuristic is the only option when very quick solutions are required, although the quality of such solutions is typically a few percentage point off the optimum. In case longer computation times are allowed, CP-SAT is the best option, able to find all the optimal solutions in less time than the RDS heuristic, that also show a small but not null optimality gap. However -- also considering the speed of the computers used for the experiments -- the exact method BIP seems the best option, since it is able to provide proven optimal solutions quickly, although the method has the advantage of taking the best of two solutions. When comparing the three exact methods, they present comparable average computation times, but the behaviour on the different problems appears different: BIP is the more robust method, with balanced times among most of the instances, although the method suffers on the larger instances; B\&B is very fast on the easier instances but it suffers heavily on the difficult instances; CP-SAT shows balanced performance but with over-long computation times on a few groups of instances that deteriorate the average performance. 

In conclusion, BIP seems the best option for all the instances apart from the largest ones, on which B\&B scales better and should be preferred. CP-SAT positions itself in between the two methods with an overall balanced behaviour.

\section{Conclusions} \label{conc}
A  formulation based on Mixed Integer Linear Programming for the Assignment Problem with Conflicts has been considered and solved via the open-source solver CP-SAT, part of the Google OR-Tools computational suite. 

The experimental results indicate that the approach we propose has performances comparable with those of the state-of-the-art exact solvers available in the literature, notwithstanding the minimal implementation effort required for our solution. Moreover, the CP-SAT solver used as a heuristic method appears to be superior to the proper heuristic methods previously introduced in the literature.

\section*{Acknowledgments}
The work was partially supported by the Google Cloud Research Credits Program.

\bibliographystyle{IEEEtran}
\bibliography{mybibfile}

\end{document}